\documentclass[11pt,english,english]{amsart}
\usepackage[latin1]{inputenc}
\usepackage{setspace}
\usepackage{amssymb}

\makeatletter

\providecommand{\LyX}{L\kern-.1667em\lower.25em\hbox{Y}\kern-.125emX\@}

 \theoremstyle{plain}    
 \newtheorem{thm}{Theorem}[section]
 \numberwithin{equation}{section} 
 \numberwithin{figure}{section} 
 \theoremstyle{remark}    
 \newtheorem*{acknowledgement*}{Acknowledgement} 
 \theoremstyle{plain}    
 \newtheorem{cor}[thm]{Corollary} 
 \theoremstyle{definition}
 \newtheorem{defn}[thm]{Definition}
 \theoremstyle{plain}    
 \newtheorem{lem}[thm]{Lemma} 
 \theoremstyle{plain}    
 \newtheorem{prop}[thm]{Proposition} 
 \theoremstyle{remark}
 \newtheorem*{rem*}{Remark}
 \theoremstyle{plain}    
 \newtheorem*{thm*}{Theorem} 

\newcommand{\hide}[1]{}

\newtheorem*{problem}{Question}

\usepackage{babel}
\makeatother
\begin{document}

\title{On processes which cannot be distinguished by finitary observation}

\thanks{This research was supported by the Israel Science Foundation (grant
No. 1333/04)}


\author{Yonatan Gutman and Michael Hochman}

\begin{abstract}
A function $J$ defined on a family $\mathcal{C}$ of stationary processes
is finitely observable if there is a sequence of functions $s_{n}$
such that $s_{n}(x_{1}\ldots x_{n})\rightarrow J(\mathcal{X})$ in
probability for every process $\mathcal{X=}(x_{n})\in \mathcal{C}$.
Recently, Ornstein and Weiss proved the striking result that if $\mathcal{C}$
is the class of aperiodic ergodic finite valued processes, then the
only finitely observable isomorphism invariant defined on $\mathcal{C}$
is entropy \cite{OW04}. We sharpen this in several ways. Our main
result is that if $\mathcal{X}\rightarrow \mathcal{Y}$ is a zero-entropy
extension of finite entropy ergodic systems and $\mathcal{C}$ is
the family of processes arising from $\mathcal{X}$ and $\mathcal{Y}$,
then every finitely observable function on $\mathcal{C}$ is constant.
This implies Ornstein and Weiss' result, and extends it to many other
families of processes, e.g. it shows that there are no nontrivial
finitely observable isomorphism invariants for processes arising from
Kronecker systems, mild and strong mixing zero entropy systems. It
also implies that any finitely observable isomorphism invariant defined
on the family of processes arising from irrational rotations must
be constant for rotations belonging to a set of full Lebesgue measure. 
\end{abstract}

\email{gyonatan@math.huji.ac.il , mhochman@math.huji.ac.il}

\maketitle
\pagestyle{myheadings}

\markboth{Y. Gutman and M. Hochman}{Classesd which cannot be distinguished by finitary invariants}

\section{Introduction}

Let $(x_{n})_{n=-\infty }^{\infty }$ be an aperiodic ergodic process
taking on finitely many values; without loss of generality the values
are in $\mathbb{N}$. We may assume that $(x_{n})$ arises from a
generating partition $\mathcal{P}=(P_{i})$ of an aperiodic, invertible
and ergodic measure preserving system $\mathcal{X}=(X,\mathcal{B},\mu ,T)$;
the system $\mathcal{X}$ is unique up to isomorphism. The question
we are interested in is: what can we learn about the underlying system
$\mathcal{X}$ by observing a sample path $(x_{n})$?

In principle, the answer is {}``everything'', since by the ergodic
theorem a typical sample path of $(x_{n})_{n=1}^{\infty }$ determines
all finite distributions of the process and this determines $\mathcal{X}$
up to isomorphism. However a more realistic scenario is one in which
at each time step another output of the process is revealed, i.e.
at time $n$ we have observed the finite sequence $x_{1}\ldots x_{n}$,
and are asked to make a guess about the nature of $\mathcal{X}$ based
on this data. 

We call a scheme for producing such a sequence of guesses an observation
scheme. To be precise,

\begin{defn}
\label{def:observation-scheme}An \emph{observation scheme} (or scheme
for short) is a metric space $\Delta $ and a sequence of functions
$s_{n}:\mathbb{N}^{n}\rightarrow \Delta $. An observation scheme
is said to \emph{converge} for a family of processes $\mathcal{C}$
if $\lim _{n\rightarrow \infty }s_{n}(x_{1}\ldots x_{n})$ exists
in probability for every process $(x_{n})\in \mathcal{C}$. A function
$J:\mathcal{C}\rightarrow \Delta $ is \emph{finitely observable}
if there is an observation scheme $(s_{n})$ which converges to $J((x_{n}))$
for every $(x_{n})\in \mathcal{C}$.
\end{defn}
Note that the larger a family of processes is, the harder it is for
a scheme to converge for every member of the family, hence large femilies
have fewer finitely observable functions.

Nonetheless, many observation schemes $(s_{n})$ are known for which
the sequence $s_{1}(x_{1})$, $s_{2}(x_{1},x_{2})$, $s_{3}(x_{1},x_{2},x_{3}),\ldots $.
converges in probability or even almost surely for \emph{every} ergodic
process $(x_{n})$. For example, if $s_{n}(x_{1}\ldots x_{n})$ counts
the frequencies of $1$'s appearing in $x_{1}\ldots x_{n}$, then
by the ergodic theorem $\lim _{n\rightarrow \infty }s_{n}(x_{1}\ldots x_{n})$
exists a.s. and equals the probability of the symbol $1$ in the process
$(x_{n})$. This example and others like it show that some things
about a process can be calculated from finite observations; but these
are generally not isomorphism invariants, and so tell us nothing about
the underlying dynamical system.

For processes $(x_{n}),(y_{n})$ etc. we denote by $\mathcal{X},\mathcal{Y}$
respectively the dynamical system determined by them. Write $(x_{n})\cong (y_{n})$
and $\mathcal{X}\cong \mathcal{Y}$ to indicate that $\mathcal{X},\mathcal{Y}$
are isomorphic as dynamical systems. We will be interested in families
of processes $\mathcal{C}$ which are closed under isomorphism, that
is, they will have the property that if $(x_{n})\in \mathcal{C}$
and $(y_{n})\cong (x_{n})$ then $(y_{n})\in \mathcal{C}$. Such a
family is called \emph{saturated}. \emph{}Usually we will specify
$\mathcal{C}$ by some property of the underlying systems, e.g. $\mathcal{C}$
might be the family of all processes arising from an irrational rotation.
In this case we would say for brevity that $\mathcal{C}$ is the class
of irrational rotations.

\begin{defn}
\label{def:isomorphism-invariants}Let $\mathcal{C}$ be a saturated
family of processes, $\Delta $ a metric space and $J:\mathcal{C}\rightarrow \Delta $.
Then $J$ is an \emph{isomorphism invariant} for $\mathcal{C}$ (or
invariant for short) if for every $(x_{n}),(y_{n})\in \mathcal{C}$,\[
(x_{n})\cong (y_{n})\; \Rightarrow \; J((x_{n}))=J((y_{n}))\]
and $J$ is a \emph{complete} invariant for $\mathcal{C}$ if the
reverse implication holds. When $J$ is an invariant we write $J(\mathcal{X})$
instead of $J((x_{n}))$. 
\end{defn}
For quite some time it has been known that the entropy $h((x_{n}))=h(\mathcal{X})$
of a process is finitely observable in the class of \emph{all} ergodic
processes. The earliest observation scheme for entropy is due to D.
Bailey \cite{Ba76}. A number of simpler schemes have been developed,
such as the Lempel-Ziv compression algorithm \cite{LZ77} and the
Ornstein-Weiss estimators \cite{OW90,OW93}.

D. Ornstein and B. Weiss recently proved a striking converse to this:
Every finitely observable invariant for the class of all ergodic processes
is a continuous function of entropy \cite{OW04}. They also showed
that there are no finitely observable invariants except entropy for
any class which contains the Bernoulli processes, for the class of
zero entropy processes or for the class of zero entropy weak mixing
processes.

However their techniques do not settle what is finitely observable
in several other interesting classes of systems. Ornstein and Weiss
have asked if there exists a complete finitely observable invariant
for the class of irrational rotations (translations by an irrational
on the group $\mathbb{R}/\mathbb{Z}$); this is not implausible, since
for this class there is a complete invariant for isomorphism, namely
the spectrum, or equivalently the modulus of rotation (up to sign
and $\bmod 1$). We remark that there are no known complete invariants
in the classes for which Ornstein and Weiss showed that entropy is
the only invariant, with the exception of the class of Bernoulli systems,
in which entropy is itself a complete invariant.

In an attempt to get a handle on this problem, we came up with the
following, which is interesting in its own right:

\begin{thm*}
Suppose $\mathcal{X}\rightarrow \mathcal{Y}$ is a zero entropy extension
of finite entropy dynamical systems, that is $h(\mathcal{X})=h(\mathcal{Y})$.
Let $\mathcal{C}$ be the class of processes arising from $\mathcal{X},\mathcal{Y}$
(that is, from generating partitions of $\mathcal{X}$ and $\mathcal{Y}$).
Then every finitely observable invariant for $\mathcal{C}$ is constant.
\end{thm*}
This allows us reclaim the results of Ornstein and Weiss, and to settle
the following problems:

\begin{thm*}
If $J$ is a finitely observable invariant on one of the following
classes:
\begin{enumerate}
\item The Kronecker systems (the class of systems with pure point spectrum) 
\item The zero entropy mild mixing processes
\item The zero entropy strong mixing processes
\end{enumerate}
Then $J$ is constant.
\end{thm*}
For the class of irrational rotations we obtain a slightly weaker
result: 

\begin{thm*}
For every finitely observable invariant $J$ on the class of irrational
rotations, there is a Borel set $\Theta \subseteq [0,1)$ of full
Lebesgue measure such that $J$ assigns the same value to processes
arising from rotations by angles in $\Theta $. In particular there
is no complete finitely observable invariant for irrational rotations.
\end{thm*}
The rest of the paper is organized as follows. Section 2 presents
some definitions and background. In section 3 we prove the theorem
about zero-entropy extensions. Section 4 contains proofs of the other
results, and in section 5 we mention some open problems.

\begin{acknowledgement*}
This paper was written as part of the authors' Ph.D. studies. We would
like to thank our advisor Professor Benjamin Weiss for his encouragement,
support and good advice.
\end{acknowledgement*}

\section{Preliminaries}

For general background on ergodic theory we refer to \cite{Hal60,Sh73,Wal82}.

\subsection{Dynamical systems, partitions and processes}

By an aperiodic ergodic system $\mathcal{X}=(X,\mathcal{B},\mu ,T)$
we mean that $(X,\mathcal{B},\mu )$ is a standard probability space,
$T$ in invertible and acts ergodically, and the set of periodic points
is of measure zero. A measure preserving systems $\mathcal{Y}=(Y,\mathcal{C},\nu ,S)$
is a \emph{factor} of the system $\mathcal{X}=(X,\mathcal{B},\mu ,T)$
if there is a measure-preserving map $f:X\rightarrow Y$ defined almost
everywhere satisfying $Sf=fT$. If there is such a map which is also
invertible and bi-measurable then $\mathcal{X}$,$\mathcal{Y}$ are
\emph{isomorphic}.

A partition $\mathcal{P}$ of $X$ is a finite ordered collection
of pairwise disjoint measurable sets $(P_{i})_{i=1}^{|\mathcal{P}|}$
whose union is $X$ (up to measure zero). If $\mathcal{P},\mathcal{Q}$
are partitions of $X$ then the partition $\mathcal{P}\lor \mathcal{Q}=(P_{i}\cap Q_{j})_{(i,j)}$
is the \emph{join} of $P,Q$ (order the pairs $(i,j)$ lexicographically);
the join of finitely many partitions is defined similarly. Write $T^{n}\mathcal{P}=(T^{n}P_{i})$. 

A partition $\mathcal{P}$ of $X$ \emph{generates} $\mathcal{X}$
if $\bigvee _{n=-\infty }^{\infty }T^{n}\mathcal{P}=\mathcal{B}$
up to measure zero, where $\bigvee _{n=-\infty }^{\infty }T^{n}\mathcal{P}$
is the $\sigma $-algebra generated by the collection $\cup _{N}\bigvee _{n=-N}^{N}T^{n}\mathcal{P}$.

For a partition $\mathcal{P}=(P_{i})_{i\in \mathbb{N}}$ and $\omega \in X$
we write $\mathcal{P}(\omega )$ for the index of the set in $\mathcal{P}$
that contains $\omega $. A partition $\mathcal{P}$ determines a
stationary ergodic process $(x_{n})$ with values in $\mathbb{N}$
by \[
x_{n}(\omega )=\mathcal{P}(T^{n}\omega )\]
 We say that $x_{i}(\omega ),x_{i+1}(\omega ),\ldots ,x_{j}(\omega )$
is the \emph{itinerary} of $\omega $ (with respect to $\mathcal{P}$)
from time $i$ to time $j$. The itinerary of $\omega $ from time
$0$ to time $N-1$ is called the $(\mathcal{P},N)$-\emph{name} of
$\omega $. If $\mathcal{P}$ is a generating partition for $\mathcal{X}$
then the system $\mathcal{X}$ and the partition $\mathcal{P}$ are
determined, up to isomorphism, by the process $(x_{n})$. We will
say this process \emph{arises} from $\mathcal{P}$ if $\mathcal{P}$
generates $\mathcal{X}$.

The space of ordered partitions of $X$ into $n$ sets comes with
a metric $\rho =\rho _{n}$ defined by \[
\rho (\mathcal{P},\mathcal{Q})=\sum _{i=1}^{n}\mu (P_{i}\triangle Q_{i})\]
for $\mathcal{P}=(P_{1},\ldots ,P_{n})$ and $\mathcal{Q}=(Q_{1},\ldots ,Q_{n})$
(here $\Delta $ denotes symmetric difference). The metric $\rho _{n}$
is complete; note however that if $\mathcal{P}_{i}\rightarrow \mathcal{P}$
in $\rho _{n}$ it may happen that some of the members of $\mathcal{P}$
are empty.

It is easy to check that if $\rho (\mathcal{P},\mathcal{Q})<\varepsilon $
then $\rho (\bigvee _{n=1}^{N}T^{n}\mathcal{P},\bigvee _{n=1}^{N}T^{n}\mathcal{Q})<N\varepsilon $.
It follows that if $\mathcal{P}_{k}\rightarrow \mathcal{P}$ in $\rho $
and $(x_{n}^{(k)}),(x_{n})$ denote the processes arising from $\mathcal{P}^{k},\mathcal{P}$
respectively, then the seuquence of processes $(x_{n}^{(k)})_{n=-\infty }^{\infty }$
converges to $(x_{n})_{n=-\infty }^{\infty }$ in probability.

Given a partition $\mathcal{P}$ of $X$ into $r$ sets and an integer
$N$ we may consider the distribution that $\mu $ induces on $\{1,\ldots ,r\}^{N}$,
where the measure of a word $w\in \{1,\ldots ,r\}^{N}$ is the measure
of the set of points whose $(\mathcal{P},N)$-name is $w$, or in
other words $\mu (\cap _{n=1}^{N}T^{-n}P_{w(n)})$. We refer to this
as the distribution of $N$-names determined by $\mathcal{P}$.

Since a distribution on $N$-names is just a $r^{N}$-dimensional
probability vector, we can compare these distributions using e.g.
the $\ell ^{1}$ metric. When we talk of closeness of $N$-name distributions,
we will mean it in this sense. Note that if $\mathcal{P},\mathcal{Q}$
are partitions and $\rho (\mathcal{P},\mathcal{Q})<\varepsilon $
then the distance between the $N$-name distributions associated with
$\mathcal{P}$ and $\mathcal{Q}$ is at most $N\varepsilon $.

\subsection{Entropy}

Let $\mathcal{X}=(X,\mathcal{B},\mu ,T)$ be an invertable ergodic
measure preserving system and $\mathcal{P}=(P_{i})$ a partition.
The \emph{entropy} of a partition $\mathcal{P}$ is\[
H(\mathcal{P})=-\sum _{i}\mu (P_{i})\log \mu (P_{i})\]
(all logarithms are to base $2$ unless specified otherwise). $H(\mathcal{P})$
is non-negative and finite (define $0\log 0=0$). The entropy of the
system $\mathcal{X}$ with respect to $\mathcal{P}$ (equivalently,
the entropy of the process arising from $\mathcal{P}$) is\begin{eqnarray*}
h(\mathcal{X},\mathcal{P}) & = & \lim _{n\rightarrow \infty }\frac{1}{n}H(\mathcal{P}\lor T\mathcal{P}\lor \ldots \lor T^{n-1}\mathcal{P})
\end{eqnarray*}
the limit above can be shown to exist. The entropy of $\mathcal{X}$
is\[
h(\mathcal{X})=\sup \{h(\mathcal{X},\mathcal{P})\, :\, \mathcal{P}\textrm{ a finite partition of }X\}\]
If $\mathcal{P}$ is a finite generating partition then $h(\mathcal{X})=h(\mathcal{X},\mathcal{P})$,
but the relation $h(\mathcal{X})=h(\mathcal{X},\mathcal{P})$ is not
in itself enough to guarantee that $\mathcal{P}$ generates. However
the Krieger generator theorem \cite{Kr70} guarentees that if $h(\mathcal{X})<\log k$
for an integer $k$ then there exists a generating partition $\mathcal{P}=(P_{1},\ldots ,P_{k})$
of $\mathcal{X}$ into $k$ sets.

In the space of partitions of $X$ into $n$ sets, the entropy is
continuous in the metric $\rho _{n}$: that is, for a partition $\mathcal{P}$,
for every $\delta >0$ there is an $\varepsilon >0$ such that if
$\rho (\mathcal{P},\mathcal{Q})<\delta $ then $|h(\mathcal{X},\mathcal{P})-h(\mathcal{X},\mathcal{Q})|<\varepsilon $.

The main fact about entropy we will use is the following classical
theorem:

\begin{thm}
(Shannon-McMillan-Breiman theorem) For any finite partition $\mathcal{P}$
of $\mathcal{X}$ and almost every $x\in X$, \[
\frac{1}{n}\log \mu (\bigcap _{i=0}^{n-1}\mathcal{P}(T^{i}x))\rightarrow h(\mathcal{X},\mathcal{P})\]

\end{thm}
A proof can be found in \cite{S96} p. 55.

Denote\[
\mu (u)=\mu (\{x\in X\, :\, \textrm{the }(\mathcal{P},n)\textrm{-name of }x\textrm{ is }u\})\]
With this notation the Shannon-McMillan-Breiman theorem states that
\[
\frac{1}{n}\log \mu (x_{1}\ldots x_{n})\rightarrow h(\mathcal{X},\mathcal{P})\]
almost surely, where $(x_{n})$ is the process arising from $\mathcal{P}$. 

Also, for partitions $\mathcal{P},\mathcal{Q}$ and $(u,v)\in \mathbb{N}^{n}\times \mathbb{N}^{n}$,
we say that $(u,v)$ is the $(\mathcal{P}\times \mathcal{Q},n)$ name
of a point $\omega \in X$ if $u$ is the $(\mathcal{P},n)$-name
of $\omega $ and $v$ is the $(\mathcal{Q},n)$-name of $\omega $.
This is just another way of talking about the partition $\mathcal{P}\vee \mathcal{Q}$.
Denote\[
\mu (v|u)=\frac{\mu (\{x\in X\, :\, \textrm{the }(\mathcal{P}\times \mathcal{Q},n)\textrm{-name of }x\textrm{ is }(u,v)\})}{\mu (\{x\in X\, :\, \textrm{the }(\mathcal{P},n)\textrm{-name of }x\textrm{ is }u\})}\]
 We will actually use the following {}``relative'' version of the
Shannon-McMillan-Breimann theorem:

\begin{thm}
\label{thm:Relative-Shannon-McMillan}(Relative Shannon-McMillan-Breimann)
Let $\mathcal{P},\mathcal{Q}$ be partitions of $\mathcal{X}$ with
entropies $h(\mathcal{X},\mathcal{P})=s\leq t=h(\mathcal{X},\mathcal{Q})$.
For every $\varepsilon >0$ there are collections of words $A_{n}\subseteq \mathbb{N}^{n}\times \mathbb{N}^{n}$
for $n=1,2,3\ldots $ such that
\begin{enumerate}
\item $\#\{u\in \mathbb{N}^{n}\, :\, (u,v)\in A_{n}\textrm{ for some }v\}<2^{(s+\varepsilon )n}$
for every $n$.
\item $\#\{v\in \mathbb{N}^{n}\, :\, (u,v)\in A_{n}\}<2^{(t-s+\varepsilon )n}$
for every $n$.
\item For almost every point $x\in X$ the $(\mathcal{P}\times \mathcal{Q},n)$-name
of $x$ is in $A_{n}$ for all sufficiently large $n$.
\end{enumerate}
\end{thm}
\begin{proof}
Define\[
A_{n}=\{(u,v)\in \mathbb{N}^{n}\times \mathbb{N}^{n}\, :\, \mu (u)>2^{-(s+\varepsilon )n}\textrm{ and }\mu (v|u)>2^{-(t-s+\varepsilon )n}\}\]
The fact that for almost every $x\in X$ the $(\mathcal{P}\times \mathcal{Q},n)$-name
of $x$ is eventually in $A_{n}$ follows from the Shannon-McMillan-Breimann
theorem, once applied to the partition $\mathcal{P}$ and once to
the partition $\mathcal{P}\times \mathcal{Q}$. The estimimates on
the size of the $u$'s represented in $A_{n}$ and the $v$'s associated
to a given $u$ in $A_{n}$ follow easily from the definition since
the mass of the $u$'s and the mass of the $v$'s relative to a given
$u$ must add to at most $1$.
\end{proof}

\subsection{Towers}

A \emph{tower of height $n$} in $\mathcal{X}$ is a set of the form
$B\cup TB\cup T^{2}B\cup \ldots \cup T^{n-1}B\subseteq X$ such that
the sets $T^{i}B$ are measurable and pairwise disjoint for $i=0,\ldots ,n-1$.
The set $B$ is called the \emph{base} of the tower, and the set $T^{i}B$
is called the $i$-th \emph{level} of the tower.

Given a partition $\mathcal{P}=(P_{i})$ and a tower $\cup _{i=0}^{n-1}T^{i}B$,
we can partition the base $B$ into disjoint (possibly empty) sets
$B_{w}$ indexed by words $w\in \mathbb{N}^{n}$, such that \[
B_{u}=\{\omega \in B\, :\, u\textrm{ is the }(\mathcal{P},n)-\textrm{name of }\omega \}\]
 This partitions the tower into disjoint subtowers $\cup _{i=0}^{n-1}T^{i}B_{u}$
whose base is $B_{u}$; these subtowers are called \emph{columns}.
Each level $T^{i}B_{u}$ is contained entirely in the element $P_{u(i)}$
of $\mathcal{P}$. Put another way, if $(x_{n})$ is the process associated
with $\mathcal{P}$ then for $\omega \in B_{u}$ the first $n$ outputs
$(x_{1}(\omega ),\ldots ,x_{n}(\omega ))$ of the process are equal
to $u=(u_{1},\ldots ,u_{n})$.

We will need two tower lemmas. 

\begin{lem}
\label{lem:infinite-tower}(Kakutani towers lemma) Let $B$ be a set
of positive measure and $N$ an integer. Then the space $X$ can be
partitioned into countably many pairwise disjoint towers all of height
no less than $N$, all of whose bases are subsets of $B$.
\end{lem}
\begin{proof}
Since $\mathcal{X}$ is aperiodic we can choose a set $B'\subseteq B$
of positive measure such that if $x\in B'$ then $T^{i}x\notin B'$
for $1\leq i<N$. Partition the base $B$ according to the first return
time to $B'$, ie let \[
B^{(n)}=\{x\in B'\, :\, n\textrm{ is the first positive integer such that }T^{n}x\in B'\}\]
Then for each $n\geq N$ we have a tower $B^{(n)}\cup TB^{(n)}\cup \ldots \cup T^{(n-1)}B^{(n)}$,
these towers are pairwise disjoint, and their union fills $X$.
\end{proof}
A stonger result is a version of the Rohlin lemmma whose proof can
be found in \cite{Sh73}

\begin{lem}
\label{lem:rohlin-tower}(Strong Rohlin lemma) Let $\mathcal{P}=\{P_{1},\ldots ,P_{k}\}$
be a partition of $X$ and $\varepsilon >0$. Then for every $N$
there is a tower $B\cup TB\cup \ldots \cup T^{N-1}B$ of height $N$
whose complement is of measure at most $\varepsilon $ and such that
the partition $\mathcal{Q}=\{B\cap P_{1},\ldots ,B\cap P_{k}\}$ induced
on $B$ by $\mathcal{P}$ has the same distribution relative to $B$
as $\mathcal{P}$ has relative to $X$.
\end{lem}
\begin{cor}
Givev $A\subseteq X$ with $\mu (A)>1-\varepsilon $ and any $N$,
there is a tower $B\cup TB\cup \ldots \cup T^{N-1}B$ in $X$ filling
all but $2\varepsilon $ of the space and with $B\subseteq A$.
\end{cor}
\begin{proof}
Let $C\cup TC\cup \ldots \cup T^{N-1}C$ be the tower provided by
the strong Rohlin lemma with respect to the partition $\{A,X\setminus A\}$
and set $B=C\cap A$.
\end{proof}

\subsection{Approximation methods for partitions}

Often a generating partition with some property is constructed by
approximation, that is, a sequence of partitions is defined satisfying
more and more of our requirements and which converge in $\rho $ to
a partition with the properties we want. Below we outline some of
the tools we use for such constructions.

If $\mathcal{A}$ is a partition or a algebra of measurable sets and
$B$ is a measurable set then we write $B\subseteq _{\varepsilon }\mathcal{A}$
to indicate that there is a set $A\in \mathcal{A}$ such that $\mu (A\triangle B)<\varepsilon $.
Clearly $B\in \mathcal{A}$ (up to measure zero) iff $B\subseteq _{\varepsilon }\mathcal{A}$
for every $\varepsilon >0$. For a partition $\mathcal{P}$ we write
$\mathcal{P}\subseteq _{\varepsilon }\mathcal{A}$ if $P_{i}\subseteq _{\varepsilon }\mathcal{A}$
for every $P_{i}\in \mathcal{P}$. 

Let $\mathcal{P}$ be a generating partition for $\mathcal{X}$ and
suppose that $\mathcal{Q}$ is a partition such that, for every $\varepsilon >0$,
there is an $N$ such that $\mathcal{P}\subseteq _{\varepsilon }\bigvee _{n=-N}^{N}T^{n}Q$.
It follows that $P\subseteq \bigvee _{n=-\infty }^{\infty }T^{n}\mathcal{Q}$,
and since $\bigvee _{n=-\infty }^{\infty }T^{n}\mathcal{Q}$ is $T$-invariant,
$\mathcal{B}=\bigvee _{n=-\infty }^{\infty }T^{n}\mathcal{P}\subseteq \bigvee _{n=-\infty }^{\infty }T^{n}\mathcal{Q}$.
Thus $\mathcal{Q}$ generates. 

Suppose $\mathcal{P},\mathcal{Q}$ are partitions of $X$ into $n$
elements and $A\subseteq _{\varepsilon }\mathcal{P}$. Then if $\rho (\mathcal{P},\mathcal{Q})<\delta $
we have $A\subseteq _{\varepsilon +\delta }\mathcal{Q}$. Thus if
$A\subseteq _{\varepsilon }\bigvee _{n=1}^{N}T^{n}\mathcal{P}$ and
$\rho (\mathcal{P},\mathcal{Q})<\delta $ then $A\subseteq _{\varepsilon +N\delta }\bigvee _{n=1}^{N}T^{n}\mathcal{Q}$. 

These observastions are essentially the proof of the following lemma,
see also \cite{Sh73} p.79:

\begin{lem}
\label{lem:approximation-principle} Let $(\mathcal{P}_{k})_{k=1}^{\infty }$
be a sequence of partitions of $X$ and $\mathcal{Q}$ a partition
of $X$. Suppose that $\rho (\mathcal{P}_{k-1},\mathcal{P}_{k})<\varepsilon (k)$
and $\mathcal{Q}\subset _{\varepsilon (k)}\bigvee _{j=-N(k)}^{N(k)}T^{-j}\mathcal{P}_{k}$
for some sequences $\varepsilon (k)>0$ and $N(k)\in \mathbb{N}$
which satisfy $\sum _{k=1}^{\infty }\varepsilon (k)<\infty $ and
$N(k)\cdot \sum _{j=k+1}^{\infty }\varepsilon (j)\rightarrow 0$ as
$k\rightarrow \infty $. Then $(\mathcal{P}_{k})$ converges to a
partition $\mathcal{P}$ and $Q\subseteq \bigvee _{j=-\infty }^{\infty }T^{-j}P$.
\end{lem}
The following theorem shows that in order to change a partition $\mathcal{P}$
into a generating partition, you need to perturb $\mathcal{P}$ by
an amount of the same order as the difference $h(\mathcal{X})-h(\mathcal{P})$.
This result is not new but we include a proof for completeness. 

\begin{thm}
\label{thm:density-of-generating-partitions}(Entropy and generating
partitions) let $h\geq 0$ and $k$ be an integer with $\log k>h$.
Let $\mathcal{X}=(X,\mathcal{B},\mu ,T)$ be an aperiodic ergodic
system with entropy $h$ and let $\mathcal{P}=(P_{1},\ldots ,P_{k})$
be a partition of $\mathcal{X}$ with $h(\mathcal{X},\mathcal{P})=h'$
(so $h'\leq h$). Then for every $\delta >0$ there is a generating
partition $\mathcal{P}'=(P'_{1},\ldots ,P'_{k})$ of $\mathcal{X}$
such that $\rho (\mathcal{P},\mathcal{P}')<\delta +\frac{h-h'}{\log k-h}$.
In particular, the generating partitions are dense in the $\rho $-metric
among the partitions of maximal entropy. 
\end{thm}
\begin{rem*}
The parameter $\delta $ was introduced only in order to deal with
the case that $h=h'$. The fact that the generating partitions are
dense among the partitions of maximal entropy is known, but we are
unable to find a reference.
\end{rem*}
\begin{proof}
Let $\delta >0$ be given. Fix a very small $\varepsilon >0$ which
will determined later. Fix a generating partition $\mathcal{Q}$ of
size $k$, and for $n=1,2,3\ldots $ let $A_{n}\subseteq \mathbb{N}^{n}\times \mathbb{N}^{n}$
be as in theorem \ref{thm:Relative-Shannon-McMillan} for the partitions
$\mathcal{P},\mathcal{Q}$ and parameter $\varepsilon $. Let $N\geq \frac{1}{\varepsilon }$
be large enough that the the set $X_{0}$ of $\omega $'s whose $(\mathcal{P}\times \mathcal{Q},n)$-name
in $A_{n}$ for all $n\geq N$ has positive measure. Applying lemma
\ref{lem:infinite-tower} we can partition the space $X$ into disjoint
towers of height at lease $\frac{N}{\varepsilon }$ whose bases are
contained in $X_{0}$, that is for each $n\geq \frac{N}{\varepsilon }$
we get disjoint towers $B^{(n)}\cup TB^{(n)}\cup \ldots \cup T^{n-1}B^{(n)}$
of height $n$ with $B^{(n)}\subseteq X_{0}$, and the union of these
towers has full measure. Partition the bases $B^{(n)}$ according
to $A_{n}$, so for a word $(u,v)\in A_{n}$ the set $B_{u,v}^{(n)}$
consists of points whose $(\mathcal{P}\times \mathcal{Q},n)$-name
is $(u,v)$. 

We construct a partition $\mathcal{P}'$ by modifying the labels of
some levels of the columns $B_{u,v}^{(n)}$. The construction proceeds
in three stages. 
\begin{description}
\item [Marking~the~base]Fix $m=\frac{1}{\varepsilon }$ (for simplicity
we ignore rounding errors and treat $m$ as an integer, and adopt
a similar philosophy later as well). Label the lower $2m$ levels
of the column $B_{u,v}^{(n)}$ (i.e. the levels indexed $0$ to $2m-1$)
with $1$'s and mark levels $2m,3m,\ldots ,[n/m]m$ with $0$'s.

The result of this procedure is that given any point $\omega \in \cup _{i=0}^{n-1}T^{i}B^{(n)}$
the base of the column can be identified as the largest index $i\in \{-n,-n+1,\ldots ,0\}$
such that the $(\mathcal{P}',2m)$-name of $T^{i}\omega $ consists
off all $1$'s. Thus given the $\mathcal{P}'$ itinerary of $\omega $
from time $-n$ to $n$, we can reconstruct the $\mathcal{P}$-name
of the column to which $\omega $ belongs. We will preserve this property
in the following steps, hence with probability $1$ given the $\mathcal{P}'$
itinerary of a point from time $-\infty $ to $-\infty $ we can determine
the $n$ corresponding to the column the point belongs to, and the
$\mathcal{P}'$-name of that column.

\item [Coding~the~$\mathcal{Q}$-itinerary~into~$\mathcal{P}'$]Denote
$A_{n}(u)=\{v\, :\, (u,v)\in A_{n}\}\subseteq \mathbb{N}^{n}$. Fix
$(u,v)\in A_{n}$ and enumerate $A_{n}(u)=\{v_{1},\ldots ,v_{r}\}$
in a way depending only on $u$; by assumption $|A_{n}(u)|<2^{(h-h'+\varepsilon )n}$.
We modify the column over $B_{u,v}^{(n)}$ so as to record the index
$i$ for which $v=v_{i}$. We do this by writing the base-$k$ representation
of $i$ near the bottom of the column. To be precise, we record the
base-$k$ digits of $i$ starting at level $2m+1$ and writing consecutively
in blocks of $m-1$, skipping levels of height $0\bmod m$ so as not
to overwrite what we did in the previous stage. Since there are at
most $2^{(h-h'+\varepsilon )n}$ possible values for $i$ we need
to overwrite $n(h-h'+\varepsilon )\log _{k}2$ levels of the column.

The result of this procedure is that if we know both the $(\mathcal{P},n)$-name
(the word $u$) and the $(\mathcal{P}',n)$-name of a point in the
base $B^{(n)}$, we can deduce its $(\mathcal{Q},n)$-name (the word
$v$) by extracting the index $i$ coded just above the base marker
in the $(\mathcal{P}',n)$ name, and looking at the $i$-th word in
the list $A_{n}(u)$. 

\item [Re-coding~the~$\mathcal{P}$-itinerary]Fix again $(u,v)\in A_{n}$.
The $\mathcal{P}$-name of the column $B_{u,v}^{(n)}$ has been partly
destroyed by the previous steps. We will fix this by overwriting still
more of the $\mathcal{P}$-name, starting where we stopped at the
previous stage, skipping levels which are at height $0\mod m$, and
stopping at some height $M=M(n)$ which we will determine. This gives
us $M-(2m+\frac{n}{m}+n(h-h'+\varepsilon )\log _{k}2)$ symbols in
which to store information. In this space we want to record the portion
of the name $u$ which has been overwritten in all three stages (including
the current stage). This consists of the first $M$ symbols of $u$
plus at most $\frac{n}{m}$ additional levels overwritten in the first
stage. Assuming as we may that $M>\varepsilon n\geq N$, we know that
the number of possibilities for the first $M$ symbols of $u$ is
bounded by $2^{(h'+\varepsilon )M}$ so using the $k$ symbols at
our disposal we need $M(h'+\varepsilon )\log _{k}2$ symbols in order
to record it, plus another $\frac{n}{m}$ symbols to record what was
erased in the first stage. Thus we require of $M$ that in addition
to $\varepsilon n<M<n$ it satisfy the inequality\[
M-(2m+\frac{n}{m}+n(h-h'+\varepsilon )\log _{k}2)\geq M(h'+\varepsilon )\log _{k}2+\frac{n}{m}\]
or equivalently\[
M\geq \frac{((h-h'+\varepsilon )\log _{k}2+2(\frac{1}{m}+\frac{m}{n}))n}{1-(h'+\varepsilon )\log _{k}2}\]
 Since $h'\leq h<\log k$, $\frac{n}{m}=\varepsilon n$ and $m=\frac{m}{n}n=\frac{1}{\varepsilon n}n\leq \frac{1}{N}n\leq \varepsilon n$,
when $\varepsilon $ is small enough it suffices that\begin{eqnarray*}
M & \geq  & \frac{((h-h'+\varepsilon )\log _{k}2+4\varepsilon )}{1-(h'+\varepsilon )\log _{k}2}n
\end{eqnarray*}
Denote the coefficiant of $n$ in expression on the right hand side
by $C(\varepsilon )$. Note that $C(\varepsilon )\rightarrow \frac{h-h'}{\log k-h'}$
as $\varepsilon \rightarrow 0$ and $0\leq C(\varepsilon )<1$. Thus
if we choose $\varepsilon >0$ small enough (in a manner depending
only on $h,h'$ and $k$) we can set $M=\max \{\varepsilon ,C(\varepsilon )\}\cdot n$
and $M$ will satisfy all the requirements, including $\varepsilon n\leq M\leq n$.

The results of this procedure is that given the $(\mathcal{P}',n)$-name
of a point in the base of the tower column $B_{u,v}^{(n)}$, we can
reconstruct its $(\mathcal{P},n)$-name by looking at the data written
in this step, and hence by the previous step its $(\mathcal{Q},n)$
name. Together with the previous stages, this means that for any point
in $X$ if we know the entire $\mathcal{P}'$ itinerary we know can
determine the column it is in and the $\mathcal{P}'$ of that column,
and hence its $\mathcal{Q}(\omega )$. This means that $\mathcal{P}'$
generates.

\end{description}
It remains to estimate how much $\mathcal{P}$ has changed. We have
modified $M+\frac{n}{m}$ levels of each column $B_{u,v}^{(n)}$,
or a $(C(\varepsilon )+\varepsilon )$-fraction of the mass of that
column. summing over all columns, this is the fraction of $X$ that
has changed. For $\varepsilon >0$ sufficiently small, this is less
than $\delta +\frac{h-h'}{\log k-h'}$, implying that $\rho (\mathcal{P},\mathcal{P}')<\delta +\frac{h-h'}{\log k-h'}$.
This completes the proof.
\end{proof}

\section{Zero-entropy extensions}

This section is dedicated to proving our main theorem, theorem \ref{thm:zero-entropy-extensions}.
Before going into the details, we would like to say a few words about
the relation of this theorem to the work of Ornstein and Weiss in
\cite{OW04}, where it was shown that entropy is the only finitely
observable invariant in some classes saturated of processes. Their
proof used a diagonalization argument: Assuming to the contrary that
for some class $\mathcal{C}$ there exists a finitely observable invariant
finer than entropy, choose two non-isomorphic processes $(x_{n}),(y_{n})\in \mathcal{C}$
with the same entropy $h$. A third process $(z_{n})$ is then constructed,
for which the observation scheme does not converge. This is done by
inductively defining the $N$-block distributions for the process
$(z_{n})$ for a sequence of rapidly increasing $N$'s, where at each
step Rohlin towers and copying lemmas are used to make $(z_{n})$
look at different time scales as though it comes from $\mathcal{X}$
or $\mathcal{Y}$. However, in order to obtain a contradiction it
must be ensured that $(z_{n})\in \mathcal{C}$, since otherwise the
observation scheme is not expected to converge. With some care one
can ensure that $(z_{n})$ is Bernoulli if $h>0$, or weak mixing
and deterministic if $h=0$, but other properties, such as pure point
spectrum or non-Bernoulliism in positive entropy, are harder to build
into $(z_{n})$.

Our results derive from the observation that when $(x_{n})$ is a
zero-entropy extension of $(y_{n})$, one can control the isomorphism
class of the diagonal process $(z_{n})$ and in fact it can be made
isomorphic to $(y_{n})$.

\begin{thm}
\label{thm:zero-entropy-extensions}Suppose $\mathcal{X}\rightarrow \mathcal{Y}$
is a zero entropy extension of finite entropy dynamical systems. Let
$\mathcal{C}$ be the family of processes arising from $\mathcal{X}$
and $\mathcal{Y}$. Then every finitely observable invariant for $\mathcal{C}$
is constant.
\end{thm}
\begin{proof}
We identify $\mathcal{Y}$ with the sub-$\sigma $-algebra of $\mathcal{X}$
which is the pull-back of the $\sigma $-algebra of $\mathcal{Y}$
through the factor map. Let $r\in \mathbb{N}$ with $\log r>h(\mathcal{X})$;
all partitions in the sequel are partitions into $r$ sets.

To simplify notation we assume that $(s_{n})$ is an observation scheme
whose range is $\mathbb{R}$; there is no loss of generality here
since given some other range we can always compose with continuous
functions from the range to $\mathbb{R}$. Suppose that there are
$\xi ,\eta \in \Delta $ such that for every pair of processes $(x_{n}),(y_{n})$
arising from $\mathcal{X},\mathcal{Y}$ respectively and generating
them, \begin{eqnarray*}
\lim s_{n}(x_{1}\ldots x_{n}) & = & \xi \quad \textrm{in probability}\\
\lim s_{n}(y_{1}\ldots y_{n}) & = & \eta \quad \textrm{in probability}
\end{eqnarray*}
We must show that $\eta =\xi $. In order to do this will construct
a generating partition $\mathcal{P}^{*}$ of $\mathcal{Y}$ and a
sequence $N(k)$ such that $s_{N(k)}(y_{1}^{*}\ldots y_{N(k)}^{*})\rightarrow \xi $
in probability (here $(y_{n}^{*})$ is the process arising from $\mathcal{P}^{*}$).
This suffices because by assumption, $\lim s_{n}(y_{1}^{*}\ldots y_{n}^{*})\rightarrow \eta $,
so $\eta =\xi $.

The partition $\mathcal{P}^{*}$ will be obtained as the limit of
a sequence of generating partitions $\mathcal{P}^{(k)}$ of $\mathcal{Y}$,
which will be constructed inductively. The induction step is provided
by the following lemma:
\begin{lem}
\label{lem:induction-step}For any generating partition $\mathcal{P}$
of $\mathcal{Y}$, and any $\varepsilon >0$, there is a generating
partition $\overline{\mathcal{P}}$ of $\mathcal{Y}$ with $\rho (\mathcal{P},\overline{P})<\varepsilon $,
and an integer $N$ so that \[
P(|s_{N}(\overline{y}_{1}\ldots \overline{y}_{N})-\xi |<\varepsilon )>1-\varepsilon \]
where $(\overline{y}_{n})$ is the process arising from $\overline{\mathcal{P}}$.
\end{lem}
Before proving the lemma let us show how it is used to prove the theorem.
We construct a sequence $\mathcal{P}^{(k)}$ of generating partitions
of $\mathcal{Y}$ and asssociated processes $(y_{n}^{(k)})$, starting
with an arbitrary generating partition $\mathcal{P}^{(0)}$ provided
by the Krieger generator theorem. 

At the induction step, given $\mathcal{P}^{(k-1)}$ we construct $\mathcal{P}^{(k)}$
using the lemma; we choose the parameter $\varepsilon =\varepsilon (k)<1/k$
in the lemma to be very small with respect to the previous stages
of the construction (see below). Thus we have\begin{equation}
\rho (\mathcal{P}^{(k-1)},\mathcal{P}^{(k)})<\varepsilon (k)\label{eq:partition-convergence-speed}\end{equation}
From the lemma we also get an integer $N(k)$ such that\begin{equation}
P(|s_{N(k)}(y_{1}^{(k)}\ldots y_{N(k)}^{(k)})-\xi |<\frac{1}{k})>1-\frac{1}{k}\label{eq:ineq-2}\end{equation}
and since $\mathcal{P}^{(k)}$ generates $\mathcal{Y}$ there is an
integer $L(k)$ such that \begin{equation}
\mathcal{P}^{(0)}\subseteq _{1/k}\bigvee _{i=-L(k)}^{L(k)}T^{i}\mathcal{P}^{(k)}\label{eq:generating-condition}\end{equation}

During the construction we are free to choose the $\varepsilon (k)$
as small as we like. First of all we will choose them so that $\sum \varepsilon (k)<\infty $.
Since the metric $\rho =\rho _{r}$ is complete (or using the Borel-Cantelli
lemma) this guarantees that $\mathcal{P}^{(k)}$ converges to a partition
$\mathcal{P}^{*}$ of $\mathcal{Y}$, with associated process $(y_{n}^{*})$.
Second, note that $\rho (\mathcal{P}^{*},\mathcal{P}^{(k-1)})\leq \sum _{m=k}^{\infty }\varepsilon (m)$.
Thus at the beginning of step $k$ of the construction, when $\mathcal{P}^{(k-1)}$
is given, we may choose a $\delta =\delta (k)>0$ depending on all
the data defined so far and prescribe that $\rho (\mathcal{P}^{*},\mathcal{P}^{(k-1)})<\delta (k)$
by requiring $\varepsilon (m)\leq 2^{-m}\delta (k)$ for every $m\geq k$.
The point is that the conditions (\ref{eq:ineq-2}) and (\ref{eq:generating-condition})
remain true for any partition (and associated process) sufficiently
close to $\mathcal{P}^{(k)}$, and hence a prudent choice of $\delta (k)$
implies that they hold for $\mathcal{P}^{*}$ and $(y_{n}^{*})$,
that is,\[
\forall m\qquad P(|s_{N(m)}(y_{1}^{*}\ldots y_{N(m)}^{*})-\xi |<\frac{1}{m})>1-\frac{1}{m}\]
and \[
\forall k\qquad \mathcal{P}^{(0)}\subseteq _{1/k}\bigvee _{i=-L(k)}^{L(k)}T^{i}\mathcal{P}^{*}\]
The first of these implies $\lim _{k\rightarrow \infty }s_{N(k)}(y_{1}^{*}\ldots y_{N(k)}^{*})=\xi $
in probability, and the second that $\mathcal{P}^{(0)}\subseteq \bigvee _{i=-\infty }^{\infty }T^{i}\mathcal{P}^{*}$,
so $\mathcal{P}^{*}$ generates $\mathcal{Y}$.
\end{proof}
{}

\begin{proof}
(of lemma \ref{lem:induction-step}) We first present a sketch of
the proof, and afterwards the details. Since $\mathcal{P}$ generates
$\mathcal{Y}$ it has full entropy, which by assumption is equal to
the entropy of $\mathcal{X}$. Therefore we can find a generating
partition $\mathcal{Q}$ for $\mathcal{X}$ with $\rho (\mathcal{P},\mathcal{Q})<\varepsilon /2$.
Let $(x_{n})$ be the process determined by $\mathcal{Q}$; then $s_{n}(x_{1}\ldots x_{n})\rightarrow \xi $
in probability, so we can choose an $N$ such that\[
P(|s_{N}(x_{1}\ldots x_{N})-\xi |<\varepsilon )>1-\varepsilon \]
Since $\mathcal{P},\mathcal{Q}$ are both defined on $\mathcal{X}$
we get a joining of the $\mathcal{P}$- and $\mathcal{Q}$-processes.
Choose now a $\delta >0$ and a suitably large $K$. Now working in
$\mathcal{Y}$ again, we can construct a partition $\mathcal{R}$
whose joint $K$-block distribution with $\mathcal{P}$ is within
$\delta $ of the joint $K$-block distribution of $\mathcal{P},\mathcal{Q}$.
Thus (assuming we chose $K$ large enough), the order of magnitude
of $\rho (\mathcal{P},\mathcal{R})$ will be of the order of $\rho (\mathcal{P},\mathcal{Q})+\delta $,
the $N$-block distribution of the $\mathcal{R}$-process will be
within $\delta $ of the $N$-block distribution of the $\mathcal{Q}$-process,
and the entropy $\mathcal{R}$ is $\delta $-close to $h(\mathcal{Y})$.
Thus although $\mathcal{R}$ doesn't necessarily generate $\mathcal{Y}$
we need only make an additional small correction to get a generating
partition $\overline{\mathcal{P}}$ for $\mathcal{Y}$, and we can
arrange that this doesn't disturb the $N$-block distributions very
much. 

Now for the details:
\begin{description}
\item [Choosing~$\mathcal{Q}$]Since $h(\mathcal{X},\mathcal{P})=h(\mathcal{Y})=h(\mathcal{X})$,
by theorem \ref{thm:density-of-generating-partitions} we can find
a generating partition $\mathcal{Q}$ for $\mathcal{X}$ with \[
\rho (\mathcal{P},\mathcal{Q})<\frac{\varepsilon }{2}\]

\item [Choosing~$N$~and~$\delta $]Denote by $(x_{n})$ the process
arising from $\mathcal{Q}$. Then $s_{n}(x_{1}\ldots x_{n})\rightarrow \xi $
in probability, so there is an integer $N$ such that \[
\mu (|s_{N}(x_{1}\ldots x_{N})-\xi |<\varepsilon )>1-\varepsilon \]
Note that condition above is a property of the $N$-block distribution
of $(x_{n})$. Thus there is a $\delta \in (0,\frac{\varepsilon }{2})$
with the property that if $(z_{n})$ is a process arising from a partition
$\mathcal{R}$ and the $N$-block distribution induced by $\mathcal{R}$
is within $\delta $ in $L^{1}$ of the $N$-block distribution of
$\mathcal{Q}$, then $\mu (|s_{N}(z_{1}\ldots z_{N})-\xi |<\varepsilon )>1-\varepsilon $.
Note also that if $\mathcal{R},\mathcal{R}'$ are two partitions of
$\mathcal{Y}$ and if $\rho (\mathcal{R},\mathcal{R}')<\delta /N$
then the $N$-bock distributions of the processes arising from $\mathcal{R},\mathcal{R}'$
differ by at most $\delta $.
\item [Choosing~$\alpha ,\beta $~and~$M$]Invoking theorem \ref{thm:density-of-generating-partitions},
choose $\alpha >0$ such that if $\mathcal{R}$ is a partition of
$\mathcal{Y}$ with entropy $h-\alpha $ then there is a \emph{generating}
partition $\mathcal{R}'$ of $\mathcal{Y}$ with $\rho (\mathcal{R},\mathcal{R}')<\delta /2N$.
Let $\beta >0$ be such that for any partition $\mathcal{S}$ of $\mathcal{Y}$,
if $\mathcal{P}\subseteq _{\beta }\mathcal{S}$ then $h(\mathcal{S})>h-\alpha $.
We may assume that $\beta <\delta /N$.

Since $\mathcal{Q}$ generates $\mathcal{X}$ and $\mathcal{P}$ is
measurable in $\mathcal{X}$ there is an $M>N$ such that \[
\mathcal{P}\subseteq _{\beta /2}\bigvee _{i=-M}^{M}T^{i}\mathcal{Q}\]
Note that this property depends only on the distribution of $(\mathcal{P}\times \mathcal{Q},2M+1)$-names,
and if $\mathcal{R}$ is a partition of $\mathcal{Y}$ such that the
distribution of $(\mathcal{P}\times \mathcal{Q},2M+1)$-names is within
$\tau $ of the distribution $(\mathcal{P}\times \mathcal{R},2M+1)$-names
(in $\ell ^{1}(\mathbb{R}^{2M+1})$) then $\mathcal{P}\subseteq _{\beta /2+\tau }\bigvee _{i=-M}^{M}T^{i}\mathcal{R}$.

\item [Choosing~$L,B$~and~$\mathcal{R}$]Fix an integer $L$ with $\max \{M,N\}/L<\beta /8$
and choose a tower $B\cup TB\cup \ldots \cup T^{L}B$ of height $L$
in $\mathcal{Y}$, filling all but $\beta /4$ of the space. We will
define a partition $\mathcal{R}$ of $\mathcal{Y}$ by modifying $\mathcal{P}$
at some of the points in the tower.

Let $(B_{u})$ be the partition of the base $B$ according to $(\mathcal{P},L)$-names.
This partition is measurable in $\mathcal{Y}$. We can further partition
each $B_{u}$ according to the $(\mathcal{Q},L)$-names as $B_{u}=\cup _{v}B_{u,v}$.
The $B_{u,v}$'s are measurable in $\mathcal{X}$ but may not be measurable
in $\mathcal{Y}$. However since $\mathcal{Y}$ is non-atomic we can
partition the sets $B_{u}$ into sets $B'_{u,v}$ in $\mathcal{Y}$
such that $\mu (B'_{u,v})=\mu (B_{u,v})$. For each $B'_{u,v}$, modify
the column over $B'_{u,v}$ so that it is labeled by $v$ (instead
of $u$). Call the resulting partition $\mathcal{R}$. 

Since \[
\rho (\mathcal{P},\mathcal{R})=2\mu (\{x\in X\, :\, \mathcal{P}(x)\neq \mathcal{R}(x)\})\]
 and on the tower $\cup _{i=0}^{L-1}T^{i}B$ we have\[
\mu \{x\in \cup _{i=0}^{L-1}T^{i}B\, :\, \mathcal{P}(x)\neq \mathcal{R}(x)\}=\mu \{x\in \cup _{i=0}^{L-1}T^{i}B\, :\, \mathcal{P}(x)\neq \mathcal{Q}(x)\}\]
and the tower fills all but $\beta /4$ of the mass, it follows that\[
\rho (\mathcal{P},\mathcal{R})\leq \rho (\mathcal{P},\mathcal{Q})+\frac{\beta }{4}<\frac{\varepsilon }{2}+\frac{\beta }{4}\]

\item [Choosing~$\overline{\mathcal{P}}$]Consider now the difference
between the distributions of $(\mathcal{P}\times \mathcal{Q},2M+1)$-names
and the distributions of $(\mathcal{P}\times \mathcal{R},2M+1)$-names.
The only difference between them is incurred at the top and bottom
$M$ levels of the tower, which have total mass $<2M/L<\beta /4$,
and the exceptional set outside the tower whose mass is $<\beta /4$.
Therefore the distributions of $(\mathcal{P}\times \mathcal{Q},2M+1)$-
and $(\mathcal{P}\times \mathcal{R},2M+1)$-names differ by at most
$\tau =\beta /2$ so\[
\mathcal{P}\subseteq _{\beta /2+\beta /2}\bigvee _{i=-M}^{M}T^{i}\mathcal{R}\]
Since the entropy of $\bigvee _{i=-M}^{M}T^{i}\mathcal{R}$ is the
same as the entropy of $\mathcal{R}$, we conclude by the choice of
$\beta $ that $\mathcal{R}$ has entropy $>h-\alpha $. We can therefore
choose a generating partition $\overline{\mathcal{P}}$ of $\mathcal{Y}$
with $\rho (\overline{\mathcal{P}},\mathcal{R})<\delta /2N$. We conclude
that\[
\rho (\mathcal{P},\overline{\mathcal{P}})<\rho (\mathcal{P},\mathcal{R})+\rho (\mathcal{R},\overline{\mathcal{P}})<\frac{\varepsilon }{2}+\frac{\beta }{4}+\frac{\delta }{2N}<\varepsilon \]
Finally, note that from the construction of $\mathcal{R}$, the $N$-block
distribution is the same as the $N$-block distribution of $Q$ except
for an error introduced by the top $N$ levels of the tower, which
have mass $<\beta /4$, and the exceptional set also of measure $\beta /4$,
which means that the $N$-block distribution of $\mathcal{R}$ and
$\mathcal{Q}$ differ by less than $\delta /2$. Since $\rho (\mathcal{R},\overline{\mathcal{P}})<\delta /2N$
we see that the $N$-block distributions of the $\mathcal{R}$-process
and the $\overline{\mathcal{P}}$-process differ by at most $\delta /2$,
so the $N$-block distributions of the $\overline{\mathcal{P}}$-process
and the $\mathcal{Q}$-process differ by at most $\delta $; by the
definition of $\delta $ this implies \[
\mu (|s_{N}(\overline{y}_{1}\ldots \overline{y}_{N})-\xi |<\varepsilon )>1-\varepsilon \]
where $(\overline{y}_{n})$ is the process defined by $\overline{\mathcal{P}}$. 
\end{description}
This completes the proof.
\end{proof}

\section{Some Applications}

An immediate consequence of theorem \ref{thm:zero-entropy-extensions}
is:

\begin{prop}
\label{pro:common-factors-and-extensions-imply-no-invariants}Let
$\mathcal{C}$ be a saturated class of processes with entropy $h$.
Suppose that every $\mathcal{X},\mathcal{Y}\in \mathcal{C}$ either
have a common factor or a common extension in $\mathcal{C}$. Then
every finitely observable invariant is constant on $\mathcal{C}$.
\end{prop}
\begin{proof}
If $\mathcal{X},\mathcal{Y}$ have a common factor $\mathcal{Z}$,
then no scheme can distinguish $\mathcal{X}$ and $\mathcal{Z}$,
and no scheme can distinguish $\mathcal{Y}$ and $\mathcal{Z}$; so
every scheme must give the same value to $\mathcal{X}$ and $\mathcal{Y}$.
The case of a common extension is similar.
\end{proof}
We turn now to some specific classes of processes. We begin by recovering
some of the results of \cite{OW04} using the techniques of the last
section.

\begin{prop}
(\cite{OW04}) There are no nontrivial finitely observable invariants
for the class of zero entropy systems or for the class of zero entropy
weakly mixing processes.
\end{prop}
\begin{proof}
Any zero-entropy ergodic systems $\mathcal{X},\mathcal{Y}$ have an
ergodic zero entropy joining (take a typical ergodic component of
$\mathcal{X}\times \mathcal{Y}$), and if $\mathcal{X},\mathcal{Y}$
are zero entropy weakly mixing systems then so is the joining $\mathcal{X}\times \mathcal{Y}$.
\end{proof}
\begin{prop}
(\cite{OW04}) If $\mathcal{C}$ is a saturated family of processes
which contains the Bernoulli processes (eg $\mathcal{C}=$all aperiodic
finite valued  ergodic processes) then entropy is the only finitely
observable invariant.
\end{prop}
\begin{proof}
For $h\geq 0$ let $\mathcal{C}_{h}=\{\mathcal{X}\in \mathcal{C}\, :\, h(\mathcal{X})=h\}$.
We must show that every finitely observable invariant scheme on $\mathcal{C}$
is constant on each $\mathcal{C}_{h}$. For $h=0$ this is the previous
proposition. For $h>0$, we use Sinai's theorem, which states that
every $\mathcal{X},\mathcal{Y}\in \mathcal{C}_{h}$ has a Bernoulli
factor with entropy $h$. By Ornsteins isomorphism theorem, these
factors are isomorphic. Since the Bernouli processes are in $\mathcal{C}$
we conclude that every $\mathcal{X},\mathcal{Y}\in \mathcal{C}_{h}$
have a common factor in $\mathcal{C}_{h}$, so every scheme is constant
on $\mathcal{C}_{h}$.
\end{proof}
Now for something new:

\begin{thm}
\label{thm:classes-without-invariants}
\begin{enumerate}
\item Every finitely observable invariant for the class of Kronecker systems
is constant
\item Every finitely observable invariant for the class of mildly mixing
zero entropy systems is constant.
\item Every finitely observable invariant for the class of strong mixing
zero entropy systems is constant.
\end{enumerate}
\end{thm}
\begin{proof}
Again, we need only note that in these classes every two systems have
a joining in the same class. 
\end{proof}
An elementary class of systems is the class $\mathcal{R}$ of irrational
rotations. A delicate and perplexing question is whether there exist
nonconstant finitely observable invariants on this class. 

To fix notation, let $([0,1),\mathcal{B},\lambda )$ be the probability
space of the unit interval with lebesgue measure. For $\alpha \in [0,1)\setminus \mathbb{Q}$
let $\mathcal{X}_{\alpha }=([0,1),\mathcal{B},\lambda ,T_{\alpha })$
where $T_{\alpha }:[0,1)\rightarrow [0,1)$ is translation by $\alpha $,
that is, $T_{\alpha }(x)=x+\alpha (\bmod 1)$. Let $\mathcal{R}=\cup \{\mathcal{X}_{\alpha }\, :\, \alpha \in [0,1)\setminus \mathbb{Q}\}$
be these systems (note that $\mathcal{X}_{\alpha }\cong \mathcal{X}_{-\alpha }$).
Thus an invariant $J:\mathcal{R}\rightarrow \Delta $ induces a map
$\widetilde{J}:[0,1)\setminus \mathbb{Q}\rightarrow \Delta $ by $\widetilde{J}(\alpha )=J(\mathcal{X}_{\alpha })$.

\begin{lem}
\label{lem:measurauility-of-rotation-invariant}If $J$ is a finitely
observable invariant on $\mathcal{R}$ then $\widetilde{J}$ is Lebesgue
measurable.
\end{lem}
\begin{proof}
We may assume that $\Delta =\mathbb{R}$ by composing continuous real-valued
functions on $s_{n}$. Let $(s_{n})$ be an observation scheme which
calculates $J$. Fix the partition $\mathcal{P}=([0,\frac{1}{2}),[\frac{1}{2},1))$
of the interval into two equal halves, and note that $\mathcal{P}$
generates for every $\mathcal{X}_{\alpha }\in \mathcal{R}$. Thus
denoting by $(x_{k}^{(\alpha )})$ the process arising from $\mathcal{P}$
and the system $\mathcal{X}_{\alpha }$, we have\[
\widetilde{J}(\alpha )=J(\mathcal{X}_{\alpha })=\lim _{n\rightarrow \infty }s_{n}(x_{1}^{(\alpha )},\ldots ,x_{n}^{(\alpha )})\]
where the limit exists in probability and is constant $\lambda $-a.e.
in $\mathcal{X}_{\alpha }$. 

Define $f_{n}:[0,1)\times [0,1)\rightarrow \Delta $ by \[
f_{n}(\alpha ,\omega )=s_{n}(x_{1}^{(\alpha )}(\omega ),\ldots ,x_{n}^{(\alpha )}(\omega ))\]
and $f:[0,1)\times [0,1)\rightarrow \Delta $ by\[
f(\alpha ,y)=\widetilde{J}(\alpha )\]

To show that $\widetilde{J}$ is measurable it suffices to show that
$f$ is measurable. And in fact, the $f_{n}$ are measurable with
respect to the product $\sigma $-algebra and since $f_{n}$ converges
in probability on every fibre $\{\alpha \}\times [0,1)$ (with respect
to $\lambda $), and the limit is the constant function $J(\alpha )$,
it follows that $f_{n}$ converges to $f$ in probability on $[0,1)\times [0,1)$
with respect to $\lambda \times \lambda $.
\end{proof}
\begin{thm}
\label{thm:class-of-rotations}Let $J:\mathbb{R}\rightarrow \Delta $
be a finitely observable invariant for $\mathcal{R}$. Then $\widetilde{J}$
is constant on a set of full measure. In particular, no finitely observable
invariant on $\mathcal{R}$ is complete.
\end{thm}
\begin{proof}
If $\alpha ,\beta \in [0,1)\setminus \mathbb{Q}$ are rationally dependent
then $\gamma =m\alpha =n\beta \in \mathbb{R}\setminus \mathbb{Q}$
for some $m,n\in \mathbb{N}$. Thus $\mathcal{R}_{\gamma }$ is a
factor both of $\mathcal{R}_{\alpha }$ and of $\mathcal{R}_{\beta }$,
so $J(\mathbb{R}_{\alpha })=J(\mathbb{R}_{\beta })$. We conclude
that $\widetilde{J}$ is a Lebesgue-measurable function on $[0,1)\setminus \mathbb{Q}$
which is constant on $\mathbb{Q}$-cosets. Any such map is constant
on a set of full measure. 
\end{proof}

\section{Remarks and probems}

Let us mention two problems which we have not been able to resolve:

\medskip{}
\noindent \begin{problem}Let $\mathcal{R}$ denote as before the
class of irrational rotations. Is every finitely observable scheme
on $\mathcal{R}$ constant? \end{problem}
\medskip{}

\medskip{}
\noindent \begin{problem}Let $\mathcal{K}$ be the class of non-Bernoulli
$K$-processes. Are there any finitely observable invariants on $\mathcal{K}$
finer than entropy?\end{problem}
\medskip{}

It has been known for some time that there are no complete Borel invariants
on $\mathcal{K}$ (the Boral structure comes from one of the natural
topologies on $\mathcal{K}$ - see Feldman's paper \cite{Feld74}).
It also follows from work of Hoffman \cite{Hoff99} that there exist
non-isomorphic $K$-systems $\mathcal{X},\mathcal{Y}$ of the same
entropy such that $\mathcal{X}\rightarrow \mathcal{Y}$ is an extension.
This implies by proposition \ref{pro:common-factors-and-extensions-imply-no-invariants}
that there are no complete finitely observable invariants on $\mathcal{K}$;
but this is not new in view of Feldman's work.

If it were true that every two processes $\mathcal{X},\mathcal{Y}\in \mathcal{K}$
had a common zero-entropy non-Bernoulli $K$-extension then proposition
\ref{pro:common-factors-and-extensions-imply-no-invariants} would
imply that there are no finitely observable invariants but entropy
on $\mathcal{K}$. However, the existence of such a joining is an
open problem.

\bibliographystyle{plain}
\bibliography{bib}

\begin{thebibliography}{10}

\bibitem{Ba76}
David Bailey.
\newblock {\em Sequential schemes for clasifying and predicting ergodic
  processes}.
\newblock Stanford University, 1976.
\newblock Ph.D. Dissertation.

\bibitem{Feld74}
Jacob Feldman.
\newblock Borel structures and invariants for measurable transformations.
\newblock {\em Proc. Amer. Math. Soc.}, 46:383--394, 1974.

\bibitem{Hal60}
Paul~R. Halmos.
\newblock {\em Lectures on ergodic theory}.
\newblock Chelsea Publishing Co., New York, 1960.

\bibitem{Hoff99}
Christopher Hoffman.
\newblock A {$K$} counterexample machine.
\newblock {\em Trans. Amer. Math. Soc.}, 351(10):4263--4280, 1999.

\bibitem{Kr70}
Wolfgang Krieger.
\newblock On entropy and generators of measure-preserving transformations.
\newblock {\em Trans. Amer. Math. Soc.}, 149:453--464, 1970.

\bibitem{OW93}
Donald Ornstein and Benjamin Weiss.
\newblock Entropy and data compression schemes.
\newblock {\em IEEE Trans. Inform. Theory}, 39(1):78--83, 1993.

\bibitem{OW04}
Donald Ornstein and Benjamin Weiss.
\newblock Entropy is the only finitely observable invariant.
\newblock {\em preprint: http://ratio.huji.ac.il/dp/dp420.pdf}, 2004.

\bibitem{OW90}
Donald~S. Ornstein and Benjamin Weiss.
\newblock How sampling reveals a process.
\newblock {\em Ann. Probab.}, 18(3):905--930, 1990.

\bibitem{Sh73}
Paul Shields.
\newblock {\em The theory of {B}ernoulli shifts}.
\newblock The University of Chicago Press, Chicago, Ill.-London, 1973.
\newblock Chicago Lectures in Mathematics.

\bibitem{S96}
Paul~C. Shields.
\newblock {\em The ergodic theory of discrete sample paths}, volume~13 of {\em
  Graduate Studies in Mathematics}.
\newblock American Mathematical Society, Providence, RI, 1996.

\bibitem{Wal82}
Peter Walters.
\newblock {\em An introduction to ergodic theory}, volume~79 of {\em Graduate
  Texts in Mathematics}.
\newblock Springer-Verlag, New York, 1982.

\bibitem{LZ77}
Jacob Ziv and Abraham Lempel.
\newblock A universal algorithm for sequential data compression.
\newblock {\em IEEE Trans. Information Theory}, IT-23(3):337--343, 1977.

\end{thebibliography}

\end{document}